\documentclass[a4paper,11pt]{article}
\usepackage{multirow}
\usepackage[english]{babel}
\usepackage[psamsfonts]{amssymb}
\usepackage{cmmib57}
\usepackage{exscale}
\usepackage{amsmath}
\usepackage{amscd}
\usepackage{latexsym}
\usepackage{graphicx}

\newtheorem{theorem}{Theorem}

\newtheorem{proposition}[theorem]{Proposition}
\newtheorem{remark}[theorem]{Remark}

\newcommand{\NN}{\mathbb{N}}

\newcommand{\PP}{\mathbb{P}}
\newcommand{\EE}{\mathbb{E}}

\newcommand{\ds}{\displaystyle}

\def\squarebox#1{\hbox to #1{\hfill\vbox to #1{\vfill}}}

\makeatother

\begin{document}

\title{On the expected number of different records in
a random sample}

\author{
Marco Ferrante\footnote{corresponding author} \ and Nadia Frigo
\\
Dipartimento di Matematica \\
Universit\`a degli Studi di Padova \\
via Trieste, 63\\
35121 Padova, Italy \\
e-mail: ferrante@math.unipd.it  \ and \ nadia.frigo@gmail.com} 
\maketitle

\begin{abstract}
Given a discrete distribution, an interesting problem is to
determine the minimum size of a random sample drawn from
this distribution,
in order to observe a given number of different records.
This problem is related with many applied problems, like the
Heaps' Law in linguistics
and the classical Coupon-collector's problem.
In this note we are able to compute theoretically
the expected size of such a sample and we provide
an approximation strategy in the case of the Mandelbrot distribution.
\end{abstract}


\section{Introduction}
Let us consider a text written in a natural language: the Heaps' law is an empirical law
which describes the portion of the vocabulary which is used in the given text.
This law can be described by the following formula
\[
R_m(n) = K n^\beta
\]
where $R_m(n)$ is the number of different words present in a text
consisting of $n$ words and
taken from a vocabulary of $m$ words, while $K$ and $\beta$
are free parameters determined empirically.
In order to obtain a formal derivation of this empirical law,
van Leijenhorst and van der Weide in \cite{Leijen_heap}
have considered the average growth in the number of records,
when elements are drawn randomly from some statistical distribution
that can assume exactly $m$ different values.
The exact computation of the average number of records
in a sample of size $n$, $\EE[R_m(n)]$, can be easily obtained
using the following approach.
Let $S=\{1,2,\ldots, m\}$ be the support of the given distribution,
define $X=m-R_m(n)$ the number of values in $S$
not observed and denote by $A_i$ the event that the
record $i$ is not observed.
It is immediate to see that $\PP[A_i]=(1-p_i)^n$,
$X=\sum_{i=1}^{m} {\bf 1}_{A_i}$
and therefore that
\begin{equation}
\label{0.5}
\EE[R_m(n)]=m-\EE[X]=m-\sum_{i=1}^{m} (1-p_i)^n \ .
\end{equation}
Assuming now
that the elements are drawn randomly from the Mandelbrot distribution,
van Leijenhorst and van der Weide
obtain that the Heaps' law is asymptotically true
as $n$ and $m$ goes to infinity and $n<<m^{\theta-1}$, where
$\theta$ is one of the parameters of the Mandelbrot distribution
(see \cite{Leijen_heap} for the details).

A slightly different problem is as follows: assume that we are
interested in the minimum number $X_m(k)$ of elements that we have to draw
randomly from a given statistical distribution in order to obtain
$k$ different records. This is clearly strictly related to the
previous problem and at first sight one expects that the technical
difficulties would be similar. However, this is not the case: in this
note we will prove that the computation of the expectation of $X_m(k)$
is more complicated and, even if related to other results
in the Coupons collector's problem, it is to the best of our knowledge
original.
The formula that we obtain is computationally hard and
we are able to perform the exact computation
in the environment R (see \cite{r-project})
just for distributions whit a support of small cardinality.
Our plan for the future is to study further this problem
in order to simplify our formula, at least in some case of interest.
By now we propose an approximation procedure in the special case
of the Mandelbrot distribution, widely used in the application, making use of the
asymptotic results proven in \cite{Leijen_heap} in order to derive the Heaps' law.

The paper is organized as follows: in the second chapter we
will derive the expected number of elements that
we have to draw from a given statistical distribution in order to
obtain $k$ different records and we will present some additional
results related to this one. In the third chapter we will compute this
value in the case of the Mandelbrot distribution. Due to the
computational effort requires to compute this expectation, we present the
exact value just for $k\le 8$. After comparing our formula
with the results obtained in \cite{Leijen_heap}, using
the exact values when $m$ is small, and the values obtained by simulation
for greater values of $m$, we use their asymptotic results
to propose an approximation to our formula.

\section{The expected value of $X_m(k)$}
Let us denote by $S=\{1, \ldots, m\}$ the
support of a given discrete distribution, by $p=(p_1,\ldots,p_m)$
its discrete density and
let us assume that the elements are drawn randomly from
this distribution in sequence.
The random variables in the sample will be independent and
the realization of each of these will be equal to $k$ with probability $p_k$.
Since we are interested here in the number of drawn one needs in order to obtain
$k$ different realization of the given distribution, let us define the following set of random variables:
$X_1$ will denote the (random) number of drawn that we need in order
to have the first record (which is trivially equal to 1),
$X_2$ will be the number of
additional drawn that we need to obtain the second
record and so on let us define, for every $i\le m$, by $X_i$ the number of
drawn needed to go from the $i-1$-th to the $i$-th different record in the sample.
From this description
we obtain that the random number $X_m(k)$
of drawn that we need to obtain $k$ different records
is equal to $X_1+\ldots + X_k$ and that $\PP[X_m(k)<+\infty]=1$.
We also define the following set of random variables:
let $Z_1$ be the type of the first record observed, $Z_2$ the type of
the second different record and so on until $Z_k$
the type of the $k$-th record observed in the sample.

\begin{remark}
\label{rm0}
The problem that we have described above is very close to the classical
Coupons collector's problem, which is usually formalized in a similar way.
In that case the random variables $X_i$ denotes the number of coupons
that we have to buy in order to go from the $i-1$-th to the $i$-th different
type of coupons in our collection and $X_m(m)$ represents the random number
of coupons that we have to collect in order to complete the collection.
The first results, due to De Moivre, Laplace and
Euler (see \cite{MR1082197} for a comprehensive introduction on this topic),
deal with the case of
constant probabilities $p_k\equiv \frac{1}{m}$,
while the first results on the unequal case have to be ascribed to Von Schelling
(see \cite{MR0061772}).
\end{remark}

In the case of a uniform distribution,
i.e. when $p_k \equiv 1/m$ for any $k\in\{1,\ldots ,m\}$,
it is immediate to see that
the random variable $X_i$, for $i\in \{2, \ldots, m\}$, has a geometric law with
parameter $(m-i)/m$. The expected number of drawn that we need in order to obtain
$k$ different records will be therefore
\begin{equation}
\label{1}
\EE[X_m(k)] = 1 + m/(m-1) +m/(m-2) + \ldots + m/(m-k+1)
\quad.
\end{equation}

When the probabilities $p_k$ are unequal,
the Coupons collector's problem fails to be useful. Indeed, in the literature it is consider
just the problem to complete the collection, i.e. in our case
to observe all the $m$ records.
This result, first proven by Von Schelling in
\cite{MR0061772}, can be obtained in
a simple and elegant way if we
look at this problem from a slightly different point of view
(see e.g. \cite{MR732623}).
Let us define the following set of random variables:
$Y_1$ will denote the (random) number of items that we need to collect to obtain the
first coupon of type $1$, $Y_2$ the number of items that we need to collect
to get the first coupon of type $2$, and so on for the others coupons.
In this setting, the waiting time to complete the collection is given by
the random variable $Y=\max\{Y_1, \ldots, Y_m\}$.
In order to compute its expected value, one
can use the Maximum-Minimums identity (see \cite{MR732623}, p.345), obtaining
\[
\begin{array}{rl}
\EE[Y] =
&
{\ds \sum_i \EE[Y_i] - \sum_{i<j} \EE[\min(Y_i,Y_j)] +
\sum_{i<j<k} \EE[\min(Y_i,Y_j,Y_k)] + \ldots}
\\
&
{\ds \ldots + (-1)^{m+1} \EE[\min(Y_1,Y_2,\ldots, Y_m)]
\quad.}
\end{array}
\]
Since the random variables $\min(Y_{i_1},Y_{i_2},\ldots, Y_{i_k})$
have a geometric law with parameter $ p_{i_1}+p_{i_2}+\ldots +p_{i_k}$, one
gets the formula
\begin{equation}
\label{2}
\EE[Y] = \sum_i \frac{1}{p_i} - \sum_{i<j}
\frac{1}{p_i+p_j} +
\sum_{i<j<k}
\frac{1}{p_i+p_j+p_k}
+ \ldots + (-1)^{m+1} \frac{1}{p_1+\ldots + p_m}
\quad.
\end{equation}

In order to compute $\EE[X_m(k)]$ for any $k\le m$, this elegant approach is
no more useful. Therefore we have to go back
to the first setting and try to compute directly the expected value of
the random variables $X_1, X_2, \ldots, X_k$.
In the case of unequal probabilities,
the law of the random variables $X_i$'s is no more so simple and,
in order to compute their expected values,
we have first to compute their conditional expected values given the
types of the preceding $i-1$-th different records obtained.
To simplify the notation, let us define
$p(i_{1},...,i_{k})=1-p_{i_1}-\ldots-p_{i_k}$ for $k \le m$
and different indexes  $i_1, i_2, \cdots , i_k$.
The main result of this section is the following proposition:
\begin{proposition}
\label{proposition1}
For any $k \in \{2, \ldots, m\}$,
the expected value of $X_k$ is equal to
\begin{equation}
\label{3}
E[X_{k}] = \sum_{i_{1}\neq i_2 \neq \cdots \neq i_{k-1}=1}^m \ \frac{p_{i_{1}}\cdots
p_{i_{k-1}}}{p(i_{1})p(i_{1},i_{2})\cdots p(i_{1},...,i_{k-1})}
\end{equation}
and therefore
\begin{equation}
\begin{array}{rl}
\label{3.5}
\EE[X_m(k)] =
&
{\ds \overset{k}{\underset{s=1}{\sum}}\EE[X_{s}] =
1+\overset{m}{\underset{i_{1}=1}{\sum}}\frac{p_{i_{1}}}{p(i_{1})}
+ \overset{m}{\underset{i_{1}\neq i_{2}=1}{\sum}}\frac{ p_{i_{1}}p_{i_{2}} }{
p(i_{1})p(i_{1},i_{2})
}+ \dots}
\\
&
{\ds
\ldots +\overset{m}{\underset{i_{1}\neq i_2 \neq \cdots \neq i_{k-1}=1}{\sum}}\frac{p_{i_{1}}\cdots
p_{i_{k-1}}}{p(i_{1})p(i_{1},i_{2})\cdots p(i_{1},\ldots ,i_{k-1})}}
\end{array}
\end{equation}
\end{proposition}
\begin{remark}
\label{rm1}
When $k=m$, expression (\ref{3.5}) represent an alternative computation
of the expected number of coupons needed to complete a collection.
The expressions (\ref{2}) and (\ref{3.5}) are different and
a direct combinatorial proof
of their equivalence seems by no means trivial.
From a computational point of view, the second formula
is heavier with respect to the first one.
In any case both of them
are not computable for large values of $k$.
\end{remark}

\noindent
\emph{Proof of Proposition \ref{proposition1}:} \
In order to compute the expected value of the variable $X_k$, we shall
conditioned this with respect to the variables $Z_1,\ldots,Z_{k-1}$,
where $Z_i$, for $i=1,\ldots,m$, denotes the type of the $i$-th different
record observed.
Let us start by evaluating $\EE[X_{2}]$:
we have immediately that
$X_{2}|Z_{1}=i$ has a (conditioned) geometric law with parameter $1-p_i=p(i)$
and therefore $\EE[X_{2}|Z_{1}=i]=1/p(i)$.
We immediately obtain that
\[
\EE[X_{2}] = \EE[\EE[X_{2}|Z_{1}]] =
\sum_{i=1}^m \ \EE[X_{2}|Z_{1}=i] \PP[Z_{1}=i]=
\sum_{i=1}^m \ \frac{p_{i}}{p(i)}
\quad.
\]
Let us now take $k \in\{3,\ldots,m\}$: it is easy to see that
\[
\EE[X_{k}]=\EE[\EE[X_{k}|Z_{1},Z_{2},\ldots ,Z_{k-1}]]
= \sum_{i_{1}\neq i_2\neq \cdots \neq i_{k-1}=1}^m
\ \EE[X_{k}|Z_{1}=i_{1},,Z_{2}=i_{2},\cdots
\]
\[
\cdots, Z_{k-1}=i_{k-1}]\ \PP[Z_{1}=i_{1},...,Z_{k-1}=i_{k-1}]
\quad.
\]
(Note that $\PP[Z_i=Z_j]=0$ for any $i\neq j$.)
The conditional law of
$X_{k}|Z_{1}=i_{1},,Z_{2}=i_{2},\cdots, Z_{k-1}=i_{k-1}$,
for $i_{1}\neq i_2\neq \cdots \neq i_{k-1}$,
is that of a geometric random variable with parameter
$p(i_{1},\ldots , i_{k-1})$ and
its conditional expected value is equal to
$p(i_{1},\ldots , i_{k-1})^{-1}$.
By the multiplication rule, we get
\[
\PP[Z_{1}=i_{1},...,Z_{k-1}=i_{k-1}]
=
\PP[Z_{1}=i_{1}]
\PP[Z_{2}=i_{2}|Z_{1}=i_{1}]\times \cdots
\]
\[
\cdots \times
\PP[Z_{k-1}=i_{k-1}|Z_{1}=i_{1},\ldots,Z_{k-2}=i_{k-2}]
\]
(note that, even though the random variables in the sample
are independent, the random variables $Z_i$
are not independent).
From its definition we have that
\[
\PP[Z_{1}=i_{1}]=p_{i_{1}}
\ ,
\]
while a simple computation gives for any $s=2,\ldots,k-1$, that
\[
\PP[Z_{s}=i_{s}|Z_{1}=i_{1},\ldots,Z_{s-1}=i_{s-1}]
= \frac{p_{i_{s}}}{1-p_{i_{1}}-\ldots -p_{i_{s-1}}}
\]
if $i_{1}\neq i_2 \neq \cdots \neq i_{k-1}$ and zero otherwise.
Recalling the compact notation $p(i_1,\ldots,i_k)=1-p_{i_1}-\ldots -p_{i_k}$,
we then get
\[
\begin{array}{rl}
E[X_{k}]
&
{\ds
= \sum^{m}_{i_{1}\neq i_2 \neq \cdots \neq i_{k-1}=1}
\frac{p_{i_{1}}p_{i_{2}}\cdots p_{i_{k-1}}}
{p(i_{1})p(i_{1},i_{2})\cdots p(i_{1},i_{2},\ldots,i_{k-1})}}
\end{array}
\]
and the proof is complete.


\begin{remark}
\label{rm2}
In the case of a uniform distribution, i.e. when $p_i\equiv 1/m$ for
any $i\in S$, we have
\[
\frac{p_{i_{1}}p_{i_{2}}\cdots p_{i_{k-1}}}
{p(i_{1})p(i_{1},i_{2})\cdots p(i_{1},i_{2},\ldots,i_{k-1})}
=
\frac{1}{(m-1)(m-2)\cdots (m-k+1)} \ .
\]
It is therefore immediate to prove that the expression (\ref{3})
coincides in this case with (\ref{1}).
\end{remark}

\section{Approximation of the expected value}
The exact formula we obtained in the
previous section is nice, but it is tremendously heavy
to compute as soon as the cardinality of the support of the distribution
becomes larger then 10.  The number of all possible ordered choices of indexes sets involved in (\ref{3.5}) increases very fast with $k$ leading to objects hard to handle with a personal computer.
For this reason it would be
important to be able to approximate this formula,
at least in some case of interest, even if its
complicated structure may suggest that it could be
quite difficult in general.
In this section we shall consider the case of
the Mandelbrot distribution, which is commonly
used in the Heaps' law and other practical problems.
Applying the results proved in \cite{Leijen_heap}, we
present here a possible strategy to approximate
the expectation of $X_m(k)$ and present some numerical
approximation in order to test our procedure.
Let us consider $R_m(n)$ and $X_m(k)$:
these two random variables are strictly related, since
$[R_m(n)>k]=[X_m(k)<n]$, for $k \le n \le m$.
However, we have seen that the computation
of their expected values is quite different.
With an abuse of notation, we could say that
the two functions $n \mapsto \EE[R_m(n)]$ and
$k \mapsto \EE[X_m(k)]$ represent one the ``inverse'' of the other.
In order to confirm this statement, let us consider the case studied in
\cite{Leijen_heap}, i.e. let us assume to sample from
the Mandelbrot distribution.
Fixed three parameters $m\in \NN$,
$\theta\in [1,2]$ and $c\ge 0$, we shall assume that
$S=\{1,\ldots,m\}$ and
\begin{equation}
\label{222}
p_i=a_m(c+i)^{-\theta}
\ \ , \ \
a_m=\left(\sum_{i=1}^{m}(c+i)^{-\theta}\right)^{-1}
\ .
\end{equation}
We implement both the expressions
(\ref{0.5}) and (\ref{3.5}) using the environment R (see
\cite{r-project}). We set the parameters of the Mandelbrot distribution
to be $c=0.30$ and $\theta=1.75$.
Using (\ref{3.5}), we compute the expected number $\mathbb{E}[X_m(k)]$
of elements we have to draw randomly from a Mandelbrot distribution in order
to obtain $k$ different records, for three levels of $m$, being $m$ the vocabulary
size, i.e the maximum size of different words. In brackets we show the expected
number of different words in a random selection of exactly $E[X_m(k)]$ elements,
computed using (\ref{0.5}). Results are collected in Table(\ref{table1}).
We see that the number of different words we expect in a text size of
dimension $E[X_m(k)]$ is close to the value of $k$ and
this supports our statement about the connection between $\EE[R_m(n)]$
and $\EE[X_m(k)]$. As underlined before, we can compute these expectations
only for small values of $k$.
\begin{table}[h]
\begin{tabular}{cc|c|c|c|l}
\cline{3-5}
& & \multicolumn{3}{c|}{Vocabulary size} \\ \cline{3-5}
& & $m=5$ & $m=8$ & $m=10$ \\ \cline{1-5}
\multicolumn{1}{|c}{\multirow{7}{*}{number of different words}} &
\multicolumn{1}{|c|}{$k=2$} & 2.80 (1.97) & 2.63 (2.00) & 2.57 (2.01)  &     \\ \cline{2-5}
\multicolumn{1}{|c}{}  &
\multicolumn{1}{|c|}{$k=3$} & 6.08 (2.87) & 5.17 (2.95)& 4.93 (2.97) &     \\ \cline{2-5}
\multicolumn{1}{|c}{} &
\multicolumn{1}{|c|}{$k=4$} & 12.42 (3.76)& 9.01 (3.90)& 8.31 (3.92)&     \\ \cline{2-5}
\multicolumn{1}{|c}{} &
\multicolumn{1}{|c|}{$k=5$} & 28.46 (4.59)& 14.81 (4.84)& 13.04 (4.88)&     \\ \cline{2-5}
\multicolumn{1}{|c}{} &
\multicolumn{1}{|c|}{$k=6$} & - & 23.95 (5.77)& 19.68 (5.84) &     \\ \cline{2-5}
\multicolumn{1}{|c}{} &
\multicolumn{1}{|c|}{$k=7$} & - & 39.96 (6.69)& 29.21 (6.80)&     \\ \cline{2-5}
\multicolumn{1}{|c}{} &
\multicolumn{1}{|c|}{$k=8$} & - & 77.77 (7.55) & 43.66 (7.74)&     \\ \cline{1-5}
\end{tabular}
\caption{Expected text size in order to have $k$ different words taken from a vocabulary of size $m$}
\label{table1}
\end{table}

\begin{figure}[h]
	\includegraphics[width=0.9\textwidth]{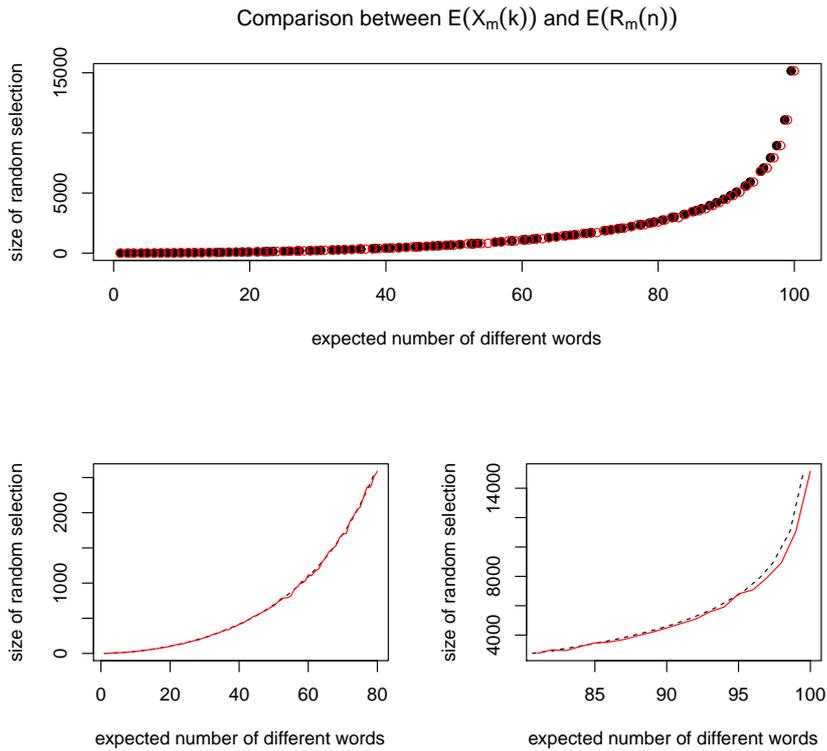}
	\caption{Comparison between $\EE[X_m(k)]$ (filled red circles) and $\EE[R_m(n)]$ (solid black circles) for $m=100$ and $k=1,\ldots,100$ (main figure). Zoom: comparison between $\EE[X_m(k)]$ (solid red line) and $\EE[R_m(n)]$ (dashed black line) for $k=1,\ldots,80$ (sx) and $k=81,\ldots,100$ (dx)}
	\label{fig1}
\end{figure}	

\begin{figure}[h]
	\includegraphics[width=0.8\textwidth]{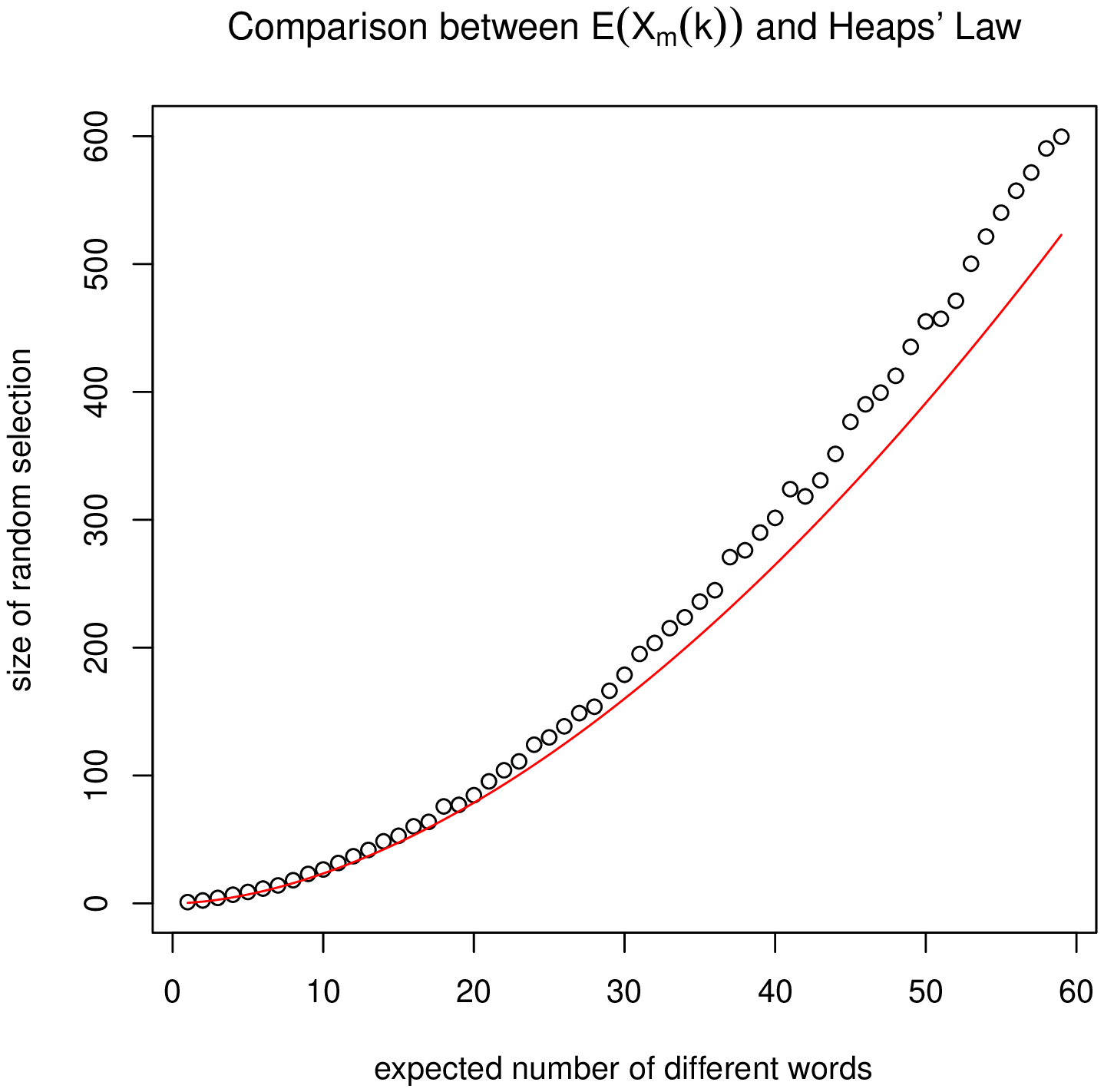}
\caption{Comparison between $\EE[X_m(k)]$ (filled black circles) and $(k/\alpha)^\theta$ (solid red line) for
	$m=500$ and $k=1,\ldots,60$}
	\label{fig2}
\end{figure}

At the same time, since $\EE[R_m(n)]\le m$,  it is clear that
our statement that $n \mapsto \EE[R_m(n)]$ and
$k \mapsto \EE[X_m(k)]$ represent one the ``inverse'' of the other
could be valid just for values of $k$ small with respect to $m$.
This idea arises also from Table (\ref{table1}), but
in order to confirm this we shall compare the two functions
for larger values of $m$. Since our formula is not
computable for values larger then $10$, we shall perform a simulation
to obtain its approximated values. In Figure (\ref{fig1}) we compare
the values of the two functions for $m=100$ and for values of $k$ ranging from $1$ to $m$.
Again, we suppose the elements are drawn from a
Mandelbrot distribution with the same value of $c$ and $\theta$.
The two functions are close up to $k=90$, while for larger values
of $k$ the distance between the two values increases.
Thanks to these results, we propose the following approximation strategy:
the main result proven in \cite{Leijen_heap} is that
\[
\EE[R_m(n)] \sim \alpha n^\beta
\]
when $n,m \rightarrow \infty$ with validity region
$n<<m^{\theta-1}$, where $\beta=\theta^{-1}$
and $\alpha=a_\infty^\beta \Gamma(1-\beta)$, where
$a_\infty=\lim_{m \rightarrow \infty} a_m$ (see the expression (\ref{222})).
Assuming that for values of $n<<m^{\theta-1}$,
$n \mapsto \EE[R_m(n)]$ and
$k \mapsto \EE[X_m(k)]$ could represent one the ``inverse'' of the other, we get
\[
\EE[X_m(k)] \sim \left(\frac{k}{\alpha}\right)^\theta
\]
with validity region
$k<<\tau$, where $\tau$ is the approximated value of $k$ for which
$\EE[X_m(k)] =m^{\theta-1}$.
In order to test our approximation scheme, we shall take
the same value of the constants as before,
$m=500$, $k=1,\ldots,60$. Figure (\ref{fig2}) shows the results: we obtain a very good correspondence
between the simulated values and the approximation curve
in the range of applicability $k<<25$.

\bibliography{coupon} 
\bibliographystyle{plain}

\end{document}